\documentclass[reqno,12pt]{amsart}

\usepackage{AKstyle}

\newcommand{\finexp}{21/22}
\newcommand{\uexp}{5/11}
\newcommand{\vexp}{6/11}
\newcommand{\hexp}{1/6}
\newcommand{\Lexp}{5/6}

\numberwithin{equation}{section}

\begin{document}

\title{A Pseduo-Twin Primes Theorem}
\author{Alex V. Kontorovich}
\email{alexk@math.brown.edu}
\address{Department of Mathematics,
Columbia University, New York, NY, $10027$} 
\curraddr{Department of Mathematics,
Brown University, Providence, RI, $02912$} 
\begin{abstract}
Selberg identified the ``parity'' barrior, that 
sieves alone cannot distinguish 
between integers having an even or odd number of factors.
We give here a short and self-contained demonstration of parity breaking using
bilinear forms, modeled on the Twin Primes Conjecture.
\end{abstract}
\date{\today}

\maketitle
\tableofcontents

\section{Introduction}
The Twin Prime Conjecture states that there are infinitely many primes $p$ such that $p+2$ is also prime.
A  refined version of this conjecture is that $\pi_2(x)$, the number of prime twins lying below a level $x$, satisfies
\[
\pi_2(x)\sim C {x \over \log ^2 x},
\]
as $x\to\infty$, where $C
 \approx 1.32032\dots$ an arithmetic constant.

The best 
result towards the Twin Prime Conjecture is Chen's 
\cite{Chen1973}, stating that there are infinitely many primes $p$ for which $p+2$ is either itself prime or the product of two primes. This statement is a quintessential exhibition of the ``parity'' barrior identified by Selberg, that sieve methods alone cannot distinguish between sets having an even or odd number of factors.  
Vinogradov's resolution \cite{Vinogradov1937}
of the ternary Goldbach problem introduced the idea that estimating certain bilinear forms can sometimes break this barrior, and there have since been many impressive instances of this phenomenon, see e.g. \cite{FriedlanderIwaniec1998, HeathBrown2001}.

In this note, we aim to illustrate 
parity breaking in a 
simple,
 self-contained  
 example. 
Consider an analogue of the Twin Prime Conjecture where 
instead of intersecting two copies of the primes, we intersect one copy of the primes with a set which analytically mimics the primes. 
For $x>2$ let 
$$
\iL(x)\sim x\log x
$$ 
denote the inverse to the logarithmic integral function, 
$$
\Li(x):=\int_2^x { dt \over \log t}
.
$$

\begin{defi}
Let $\hat{\pi}(x)$ denote the number of primes $p\le x$ such that $p=\lfloor \iL(n) \rfloor$
for some integer $n$.  
\end{defi}

Here $\lfloor \cdot\rfloor$ is the floor function, returning the largest integer not exceeding its argument.
Our main goal is to demonstrate
\begin{theorem}\label{thm:main}
As $x\to\infty$,
\[
\hat{\pi}(x) \sim  \frac{x}{\log^2 x}.
\]
\end{theorem}
Notice that the constant above is $1$, that is, there is no arithmetic interference.
This theorem follows also from the work of Leitmann \cite{Leitmann1977}; both his proof and ours essentially mimic  Piatetski-Shapiro's theorem \cite{PiatetskiShapiro1953}. Our aim is to 
give a short derivation of this statement 
 from scratch.

\subsection*{Outline}
In \S\ref{sec:Ests} we give 
bounds for exponential sums of 
linear and bilinear type; 
these are used
in the sequel.
We devote
\S\ref{sec:Reduct}
to
reducing Theorem \ref{thm:main} to 
an estimate for exponential sums over 
primes. 
The latter are treated in
\S\ref{sec:Vaughan} by Vaughan's identity,
relying on
the bounds 
of 
\S\ref{sec:Ests}
to establish Theorem \ref{thm:main}.

\subsection*{Acknowledgements}
The author wishes to thank Peter Sarnak for suggesting this problem, and generously lending of his time.
Thanks also to Dorian
Goldfeld and Patrick Gallagher for enlightening conversations, and to
Tim Browning, 
Gautam Chinta,
 Steven J. Miller, 
and the 
referee for helpful remarks and corrections to an earlier draft.

\section{Estimates for Linear and Bilinear Sums}\label{sec:Ests}

In this section we develop preliminary bounds of linear and bilinear type, which are used in the sequel.
We require first the following two well-known estimates due originally to Weyl \cite{Weyl1921}  and van der Corput \cite{vanderCorput1921,vanderCorput1922}; see e.g. Theorem 2.2 and Lemma 2.5 of
\cite{GrahamKolesnik1991}.

\begin{lemma}[\protect van der Corput] \label{thm:vdC}
Suppose $f$ has two continuous derivatives and for $0<c<C$, we have $c\Delta \le f'' \le C
\Delta$ on $[N,2N]$. Then \benn \sum_{N<n\le N_1 \le 2N} e(f(n))\ \  \ll_{C,c} \ \ 
N \Delta^{1/2} + \Delta^{-1/2}.\eenn
\end{lemma}

This is proved by truncating Poisson summation, comparing the sum to the integral, and integrating by parts two times.

\begin{lemma}[Weyl, van der Corput]\label{WvdC}
Let $z_k\in\C$ be any complex numbers, $k=K+1,\dots,2K$. Then for any $Q\le K$,
\benn
\left|\sum_{K<k\le2K} z_k \right|^2 \le \frac{K+Q}{Q}\sum_{|q|<Q}
(1-\frac{|q|}{Q}) \sum_{K<k,k+q\le2K}z_k \bar{z}_{k+q} 
.
\eenn 
\end{lemma}

To prove this, shift the interval by $q$ and average the contributions over $|q|<Q$.

\subsection{Estimating Type I Sums}

We use Lemma \ref{thm:vdC} to prove 

\begin{lemma}\label{lemma:exp}
For any integer $h$ and $\ell\ge 1$,
\benn 
\sum_{N< n\le N_1 \le 2 N} e(h\Li(n \ell)) \ll \twocase{}{N}{if
$h=0$}{ (N|h|\ell)^{\foh} \log(N \ell)}{otherwise.} 
\eenn
\end{lemma}
\begin{rmk}
Here as throughout, the implied constant is absolute unless otherwise specified.
\end{rmk}
\begin{proof}
Let $\Xi$ denote the sum in question. The trivial estimate is $
N$. Assume without loss of generality $h> 0$. Apply 
Lemma \ref{thm:vdC}
with $f(n)=h\Li(n \ell)$, taking $\Delta = \frac{h
\ell}{N \log^2 (N \ell)}$. Thus \benn \Xi \ll N \left(\frac{h \ell}{N \log^2 (N
\ell)}\right)^{1/2} + \left(\frac{N \log^2 (N \ell)}{h \ell} \right)^{1/2} \ll (N h \ell)^{1/2}
\log (N \ell),\eenn so we are done.
\end{proof}

\subsection{Estimating Type II Sums}

We first require the following estimate.

\begin{lemma}\label{lemma:So}
For positive integers 
$h,k, q,$ and  $L\ge10$, 
let 
\be\label{eq:So}
S_0(q;k) := \sum_{L<\ell\le 2L} e\bigg(h\big[\Li(\ell k)-\Li(\ell(k+q))\big]\bigg) 
. 
\ee 
Then
\benn 
S_0(q;k) \ll (L h q)^{1/2} 
,
\eenn
where the implied constant is absolute, that is, independent of $k$.
\end{lemma}
\begin{proof}

We again apply Lemma \ref{thm:vdC}
, this time choosing for the function $f(x) = h(\Li(x
k) - \Li(x(k+q)))$. Then for $L<x\le 2L$, 
\benn 
f''(x) = h\left( \frac{-k}{x \log^2(x
k)} +\frac{k+q}{x \log^2(x (k+q))} \right) 
= 
h q \frac{\log(k'
x)-2}{x \log(k' x)}
\asymp
{hq \over L}
,
\eenn 
for some $k'\in [k,k+q)$ by the Mean
Value Theorem in $k$. Thus we can take $\Delta = \frac{h q}{L}$ and
\benn S_0 \ll L \left(\frac{h q }{L}\right)^{1/2}+\left(\frac{L}{h q}\right)^{1/2} \ll (L
h q)^{1/2}, \eenn as desired.
\end{proof}

With this estimate in hand, we control Type II sums as follows (see also \cite[Lemma 4.13]{GrahamKolesnik1991}).

\begin{lemma}\label{lemma:II}
Let $\ga(\ell)$ and $\gb(k)$ be sequences of complex numbers supported
in $(L,2L]$ and $(K,2K]$, respectively, and suppose that 
\be\label{eq:L2norms}
\sum_\ell
|\ga(\ell)|^2 \ll L\log^{2A} L
\quad\text{ and }\quad
\sum_k |\gb(k)|^2 \ll K\log^{2B} K
.
\ee
Then 
\be\label{eq:Lem2} 
\sum_{L<\ell\le 2L}\sum_{ K<k\le 2K} \ga(\ell) \gb(k) e(h\Li(\ell
k)) \ll K L^{\Lexp}  h^{\hexp} \log^{A}L \ \log^{B}K 
. 
\ee
The implied constant in \eqref{eq:Lem2} depends only on those in \eqref{eq:L2norms}. 
\end{lemma}
\begin{proof}
Let $S$ denote the sum on the left hand side. By Cauchy-Schwartz,
\benn |S|^2 \ll \left(\sum_\ell |\ga(\ell)|^2\right) \sum_\ell \left| \sum_k \gb(k)
e(h\Li(\ell k))  \right|^2. \eenn

Let $Q\le K$ be a parameter to be chosen later. 
%
%
%
Using 
%
Lemma \ref{WvdC}
and 
\eqref{eq:L2norms},
we get:
\beann
|S|^2 &\ll& L \log^{2A} L \ \frac{K+Q}{Q} \sum_{|q|<Q} (1-\frac{|q|}{Q}) \\
& & \times\sum_\ell \sum_{K<k, k+q\le 2K} \gb(k) \bar{\gb}(k+q)\
e\bigg(h\big[\Li(\ell k)-\Li(\ell(k+q))\big]\bigg)\\
 & \ll & L \log^{2A} L \  \frac{K}{Q} \sum_{1\le |q| <Q} \sum_k |\gb(k) \bar{\gb}(k+q)| |S_0(q;k)| + \frac{K^2L^2}{Q}\log^{2A} L \ \log^{2B} K,
\eeann where $S_0$ is defined by \eqref{eq:So}.

Using Cauchy's inequality, $|x\bar{y}|\le \foh (|x|^2 +
|y|^2)$, and the fact that $|S_0(q;k)| =| S_0(-q;k+q)|$, we get

\beann |S|^2 &\ll&  \frac{K^2L^2}{Q}\log^{2A}L \ \log^{2B}K + \frac{L
K}{Q} \log^{2A} L \sum_k |\gb(k)|^2 \sum_{1\le q <Q} |S_0(q;k)|. \eeann
From Lemma \ref{lemma:So} we have the estimate: \benn \frac{1}{Q}
\sum_{1\le q <Q} |S_0(q;k)| \ll (L h Q)^{1/2}, \eenn so we finally
see that \benn |S|^2 \ll \frac{K^2L^2}{Q}\log^{2A}L \ \log^{2B}K + L^{3/2}
K^2 \log^{2A} L \ \log^{2B}K  Q^{1/2} h^{1/2}. \eenn The choice $Q=\lfloor
L^{1/3}h^{-1/3}\rfloor$ gives the desired result.
\end{proof}

\section{Reduction to Exponential Sums}\label{sec:Reduct}

In this section, we reduce the statement of Theorem \ref{thm:main} to a certain exponential sum over primes.
We follow standard methods, see e.g. \cite{GrahamKolesnik1991, HeathBrown1983}, which we include here for completeness.
If $p=\lfloor\iL(n)\rfloor$ then $p\le \iL(n) < p+1$, or
equivalently, $\Li(p) \le n < \Li(p+1)$. The existence of an integer
in the interval $[\Li(p),\Li(p+1))$ is indicated by the value
$\lfloor\Li(p+1)\rfloor-\lfloor\Li(p)\rfloor$, so we have 
\benn
\hat\pi(x) = \sum_{p\le x} 
\bigg(
\lfloor\Li(p+1)\rfloor-\lfloor\Li(p)\rfloor 
\bigg)
. 
\eenn 
Write
$\lfloor\gt\rfloor=\gt - \psi(\gt)-\foh$, where $\psi$ is the shifted
fractional part 
\[
\psi(\gt):=\{\gt\}-\foh \in [-\foh,\foh).
\] 

So we have: 
\benn 
\hat\pi(x) 
= 
\sum_{p\le x} 
\bigg[\Li(p+1)-\Li(p)\bigg] 
+
\sum_{p\le x} 
\bigg[
\psi(\Li(p))-\psi(\Li(p+1))
\bigg]
. 
\eenn 
Since $\Li
'(x)=\frac{1}{\log x}$, we use the Taylor expansion: 
\[
\Li(p+1) =
\Li(p) + \frac{1}{\log p} + O\left(\frac{1}{p\log^2p}\right)
\]
 to get: 
\benn
\hat\pi(x) = \sum_{p\le x} \frac{1}{\log p} + \sum_{p\le x} \bigg[ \psi(
\Li(p)) -\psi( \Li(p+1)) \bigg]+ O(1). 
\eenn 
By partial summation and a crude form of the 
Prime Number Theorem, 
\beann 
\sum_{p\le x} \frac{1}{\log p}& =&
\int_2^x \frac{d\pi(t)}{\log t} = \frac{\pi(x)}{\log x} + O
\left(\int_2^x \frac{\pi(t)}{t \log^2 t} dt\right)  \\
& =& \frac{x}{\log^2x}+
O\left(\frac{x}{\log^3 x}\right).
\eeann 
Therefore to prove Theorem \ref{thm:main}, it suffices to show that 
\be\label{reduce1}
\sum_{p\le x} 
\bigg[
\psi( \Li(p)) -\psi( \Li(p+1))
\bigg]
 \ll \frac{x}{\log^3 x}.
\ee 
Equivalently,  split  the sum into dyadic segments and apply partial summation to reduce \eqref{reduce1} to the statement that for any $N<N_{1}\le2N$,
\be\label{eq:gSDef1}
\Sigma := \sum_{N<n\le
N_1\le 2N} 
\gL(n)
\bigg[
\psi(\Li(n))-\psi(\Li(n+1))
\bigg] 
\ll
\frac{N}{\log^2 N},
 \ee
with $N\ll x$.
Here $\gL$ is the von Mangoldt function:
$$
\gL(n)=\twocase{}{\log p}{if $n=p^k$ is a prime power}{0}{otherwise.}
$$

The truncated Fourier series of $\psi$ is
 \be\label{eq:fracExpand} 
 \psi(\gt)=\sum_{0<|h|\le
H} c_h\ e(\gt h) + O(g(\gt,H)),
 \ee 
where 
$e(x)=e^{2 \pi i x}$,
$c_h = \frac{1}{2\pi ih}$, and
$$
g(\gt,H)=\min\left(1,\frac1{H\|\gt\|}\right)
.
$$
 Here $\|\cdot\|$ is the distance to the nearest integer. 
 In the above, $H$ is a parameter which we will
choose later, eventually setting 
\be\label{eq:His}
H=\log^4N.
\ee 
The function $g$ has Fourier expansion 
$$
g(\gt,H)
=
\sum_{h\in\Z}
a_{h}\
e(\gt h)
,
$$
in which
\be\label{eq:ahBnd}
a_{h}
\ll
\min\left(\frac{\log2H}H,\frac H{|h|^{2}}\right)
.
\ee

Using \eqref{eq:fracExpand}  write the
sum in \eqref{eq:gSDef1} as $\Sigma=\Sigma_1 + O(\Sigma_2)$ where 
\benn
 \Sigma_1 :=
\sum_n \gL(n) \sum_{0<|h|\le H}  c_h
\bigg[
e(h \Li(n)) -
e(h\Li(n+1))
\bigg] 
\eenn 
and 
\benn 
\Sigma_2 :=
\sum_{n\sim N}\gL(n)\bigg[
g(\Li(n),H)+g(\Li(n+1),H)
\bigg] . 
\eenn

We first dispose of $\gS_{2}$. Using positivity of $g$, the bound \eqref{eq:ahBnd}, and Lemma \ref{lemma:exp}, we have
\beann
\gS_{2}
&\ll&
\log N\sum_{n\sim N}g(\Li(n),H)
\ll 
\log N
\sum_{h\in\Z}
|a_{h}|
\left|
\sum_{n\sim N}
e(\Li(n)h)
\right|
\\
&\ll&
\log N
\left[
{\log 2H\over H}
N
+
\sum_{h\neq 0}
{H\over |h|^{2}}
(N|h|)^{1/2}
\log N
\right]
\\
&\ll&
(\log N)^{2}
\bigg(
N/H
+
N^{1/2}
H
\bigg).
\eeann
This bound is acceptable for \eqref{eq:gSDef1} on setting $H$ according to \eqref{eq:His}.

Next we massage $\gS_{1}$.
On writing $\phi_h(x)= 1 - e(h(\Li(x+1)-\Li(x))$, we  see by partial
summation  that
\beann
\Sigma_1 & \ll & \sum_{1\le h \le H} h^{-1} \left| \sum_{N<n\le N_1} \gL(n) \phi_h(n)e(h\Li(n)) \right| \\
 & \ll & \sum_{1\le h \le H} h^{-1}  \left| \phi_h(N_1)  \sum_{N<n\le N_1} \gL(n)e(h\Li(n))\right| \\
 & & + \int_N^{N_1}  \sum_{1\le h \le H} h^{-1}\left| \frac{\partial \phi_h (x)}{\partial x} \sum_{N<n\le x}  \gL(n)e(h\Li(n))\right| dx
\\
 & \ll & \frac{1}{\log N} \max_{N_2 \le 2N} \sum_{1\le h \le H} \left|\sum_{N<n\le N_2}\gL(n)e(h\Li(n))\right|.
\eeann
Here we used the bounds
\benn
\phi_h(x) \ll
h(\Li(x+1)-\Li(x)) \ll \frac{h}{\log N}
\eenn
and
\benn
\frac{\partial \phi_h (x)}{\partial x} \ll h \left(\frac{1}{\log(x+1)} -
\frac{1}{\log(x)}\right)\ll \frac{h}{N \log^2 N}
\eenn
for $N\le x \le
2N$. We have thus reduced Theorem \ref{thm:main} to the statement that for all $N<N_{2}\le2N$,
\be\label{eq:left}
S:=\sum_{0<h\le H} \left| \sum_{N<n\le N_2\le2N}
\gL(n) e(h\Li(n))  \right| \ll \frac{N}{\log N}.
\ee
We establish this fact
in the next section.


\section{Proof of Theorem \ref{thm:main}
}\label{sec:Vaughan}

Our goal in this section is to demonstrate \eqref{eq:left}, thereby establishing Theorem \ref{thm:main}. We will actually prove more; instead of a log savings, we will save a power:

\begin{theorem}
For $S$ defined in \eqref{eq:left} and any $\gep>0$, we have 
\[
S \ll_{\gep} N^{\finexp+\gep}, \qquad\text{ as }\qquad N\to\infty
.
\]
\end{theorem}

Fix $u$ and $v$,
parameters to be chosen later, and let $F(s)=\sum_{1\le n \le v}
\gL(n) n^{-s}$ and $M(s)=\sum_{1\le n \le u} \mu(n) n^{-s}$, where $\mu$ is the M\"obius function:
$$
\mu(n)=\twocase{}{(-1)^k}{if $n$ is the product of $k$ distinct primes}{0}{if $n$ is not square-free.}
$$ 
The functions $F$ and $M$ are the truncated Dirichlet polynomials of the functions $-\gz'/\gz$ and $1/\gz$, respectively, where $\gz(s)$ is the Riemann zeta function.  Notice,
for instance, that 
\[
\frac{\zeta'}{\zeta}(s)+F(s) = -\sum_{n>v} \gL(n)
n^{-s}. 
\]
Comparing the Dirichlet coefficients on both sides of the
identity
\benn \frac{\zeta'}{\zeta}+F = \left(\frac{\zeta'}{\zeta}+F
\right)(1-\zeta M) + \zeta'M + \zeta FM \eenn 
gives for $n>v$:
\be\label{eq:Vaughan}
- \gL(n) = - \sum_{{k\ell=n}\atop{k>v, \ell>u}}\gL(k)
\sum_{{d|\ell}\atop{d>u}} \mu(d) - \sum_{{k\ell=n}\atop{\ell\le u}}\log k \
\mu(\ell) + \sum_{{k\ell m=n}\atop{\ell\le v, m\le u}}1\cdot \gL(\ell) \mu(m)
\ee
This formula is originally due to Vaughan \cite{Vaughan1977} (see also \cite[Lemma 4.12]{GrahamKolesnik1991}).
Assume for now that $v\le N$ (we will eventually set $u$ and $v$ to
be slightly less than $\sqrt{N}$). Multiply the above identity by
$e(h\Li(n))$ and sum over $n$: 
\beann
\sum_{N<n\le N_2\le 2N} \gL(n) e(h\Li(n)) 
&=& \sum_{u<\ell\le N_2 /v} \sum_{{N/\ell\le k \le N_2/\ell}\atop{v<k}} \gL(k) a(\ell) e(h\Li(k\ell)) \\
& & + \sum_{\ell\le u} \sum_{N/\ell \le k\le N_2/\ell} \mu(\ell)\log k \ e(h\Li(k\ell)) \\
& & - \sum_{r\le uv} \sum_{N/r \le k \le N_2/r} b(r) e(h \Li(k r))\\
&=& S_1 +S_2 - S_3, \eeann where \benn a(\ell)= \sum_{{d|\ell}\atop{d>u}}
\mu(d) \text{, and } b(r)=\sum_{{\ell m=r}\atop{\ell\le v, m\le u}}\gL(\ell)
\mu(m). \eenn

It is the bilinear nature of the above identity which we exploit, forgetting the arithmetic nature of the coefficients $a$, $b$, $\mu$, and $\gL$, and just treating them as arbitrary. The savings then comes from the matrix norm of $\{e(h\Li(kl))\}_{k,\ell}$. This is achieved as follows.

Notice that $|a(\ell)|$ is at most $d(\ell)$, the number of divisors of $\ell$, and similarly $|b(r)|\le
\sum_{d|r}\gL(d) = \log r$, so we have the estimates

\benn \sum_{L<\ell\le 2L}|a(\ell)|^2 \ll L \log^3 L \text{, and } \sum_{R
< r \le 2R}|b(r)|^2 \ll R\log^2 R. \eenn

It now suffices to show that $\sum_{0<h<H} |S_i| \ll_{\gep}
N^{\finexp+\gep}$ for each $i=1,2,3$ by choosing $u$ and $v$
appropriately. We treat the sums of $S_i$ individually in the next
three subsections.

\subsection{The sum $S_2$}

Let $G(x):=\sum_{k\le x} e(h\Li(k\ell))$. By Lemma \ref{lemma:exp},
$G(x) \ll (x h \ell)^{\foh}\log(x\ell)$, so by partial integration we get
\beann
S_2 &=&  \sum_{\ell\le u}\mu(\ell) \sum_{N/\ell \le k\le N_2/\ell} \log k \ e(h\Li(k\ell)) \\
& \ll &  \sum_{\ell\le u} \left| \int_{N/\ell}^{N_2/\ell} \log x \ dG(x) \right| \\
& \ll & \sum_{\ell\le u} \left( \sqrt{N h}\log^2 N + \int_{N/\ell}^{N_2/\ell} \frac{1}{x}\sqrt{x h \ell}\log (x \ell) dx\right)\\
& \ll & \sqrt{N h } u \log^2 N. 
\eeann 
Thus $\sum_{1\le h <H} |S_2|
\ll_{\gep} N^{\finexp+\gep}$ (as desired) on taking $u=N^{\uexp}$ and
recalling \eqref{eq:His}. 

\subsection{The sum $S_1$}

Rewrite $S_1$ and split it into $\ll \log^2N$ sums of the form:
\beann
S_1 & = & \sum_{{N\le k\ell \le N_2}\atop{v<k, u<\ell}} \ga(k) \gb(\ell) e(h\Li(k\ell)) \\
& \ll & \log^2N \sum_{L<\ell\le 2L}\sum_{{K<k\le 2K}\atop{N<k\ell\le N_2}}
\ga(k) \gb(\ell) e(h\Li(k\ell)). \eeann

The roles of $k$ and $\ell$ are essentially symmetric (allowing $\ga$
and $\gb$ to be either $\gL$ or $a$ affects only powers of $\log$
and not the final estimate) and taking $v=u$, we may arrange it so
$N^{\uexp} \le K \le N^{1/2} \le L \le N^{\vexp}$.

Now using Lemma \ref{lemma:II}, we find that:
\beann S_1 &\ll& \log^2N \left( K L^{\Lexp} h^{\hexp} \log^2L \
\log^2K \right)
\\
& \ll & \log^6N \left( N^{\finexp} h^{\hexp} \right). \eeann Thus
$\sum_h |S_1| \ll_{\gep} N^{\finexp+\gep}$ as desired.

\subsection{The sum $S_3$}

Recall $S_3$ and break it according to: \beann
S_3 &= & \sum_{r\le uv} b(r)  \sum_{N/r \le k \le N_2/r} e(h \Li(k r))\\
 & = & \sum_{r\le u} + \sum_{u<r\le uv}\\
& = & S_4 + S_5. \eeann

We treat $S_4$ exactly as $S_2$, getting $S_4\ll (N h)^{1/2}  \log N
(u \log u)$, which is clearly sufficiently small.

For $S_5$, the analysis is identical to that of $S_1$ and gives the
same estimate, so we are done.


\bibliographystyle{alpha}

\bibliography{AKbibliog}

\end{document}